\documentclass{article}

\usepackage[english]{babel}
\usepackage{amsmath}
\usepackage{amsfonts}
\usepackage{amssymb}
\usepackage{amstext}
\usepackage{setspace}
\usepackage{graphicx}
\usepackage{bm}
\usepackage{mathrsfs}
\usepackage{float}
\usepackage{authblk}
\usepackage{esint}

\begin{document}

\title{Towards a proof of the Riemann Hypothesis \\ Searching for collaborators}
\author[1]{Guilherme Rocha de Rezende}
\affil[1]{Federal Institute of Brasília-Brazil\footnote{E-mail: {guilherme.rezende@ifb.edu.br}}}

\date{\today}

\maketitle

\begin{abstract}
\begin{center}
In this article we propose a revisitation of the well-known argument principle that may lead to the solution of the Riemann hypothesis. We are looking for collaborators.
\end{center}
\end{abstract}

\section{Introduction}
 The Riemann Hypothesis is a famous conjecture made by Bernhard Riemann in his article on prime numbers. Riemann, as indicated by the title of his article \cite{Riemann}, wanted to know the number of prime numbers in a given interval of the real line, so he extended a Euler observation and defined a complex function called Riemann zeta function. Riemann obtained an explicit formula, which depends on the non-trivial zeros of the zeta function, for the quantity he was looking for. Along the way, Riemann mentions that probably all non-trivial zeros of the zeta function are, in the now called critical line, that is, when the complex argument $ s = \sigma + IT $ of the zeta function has a real part equal to one-half. - $ \sigma = \frac{1}{2}$.
 
Despite being one of the most famous and extensively studied problems in mathematics, the Riemann Hypothesis remains unsolved. However, progress has been made in understanding its properties and implications, and many new techniques and ideas have emerged from the study of this important problem. In fact, the Riemann Hypothesis has been the driving force behind many significant developments in mathematics over the past century and a half, and it continues to be a central focus of research in the field of number theory today. In this article we propose a methodology that will possibly lead to the solution of the Riemann hypothesis.

\newpage

\section{Reasoning}
Backlund's method \cite{Blacklund} is based on Riemann's observation that if $N(t)$ denotes the number of roots, of xi function $\xi(s)$, $\rho$ in the range$\{ 0 < \Im s < t\}$, then \footnote{This follows from the argument principle of complex analysis or, more directly, from termwise integration of the uniformily convergent series $\frac{\xi'(s)}{\xi(s)}=\sum(s-\rho)^{-1}$ using the  Cauchy integral formula.}
\begin{align}
N(t)=\frac{1}{2 \pi i} \ointctrclockwise_{\partial R} \frac{\xi'(s)}{\xi(s)}ds \label{1}
\end{align}
where $\xi(s)$ is
\begin{align}
\xi(s)=\frac{1}{2}s(s-1)\pi^{-s/2}\Gamma(s/2)\zeta(s).\label{2}
\end{align}
$R$ is a rectangle of the form $\{-\epsilon \leq \Re s \leq 1 + \epsilon, 0 \leq \Im s \leq t\}$, and $\partial R$ is the boundary of R oriented in the usual counterclockwise direction,
\begin{figure}[h]
\centering
\includegraphics[width=0.6\textwidth]{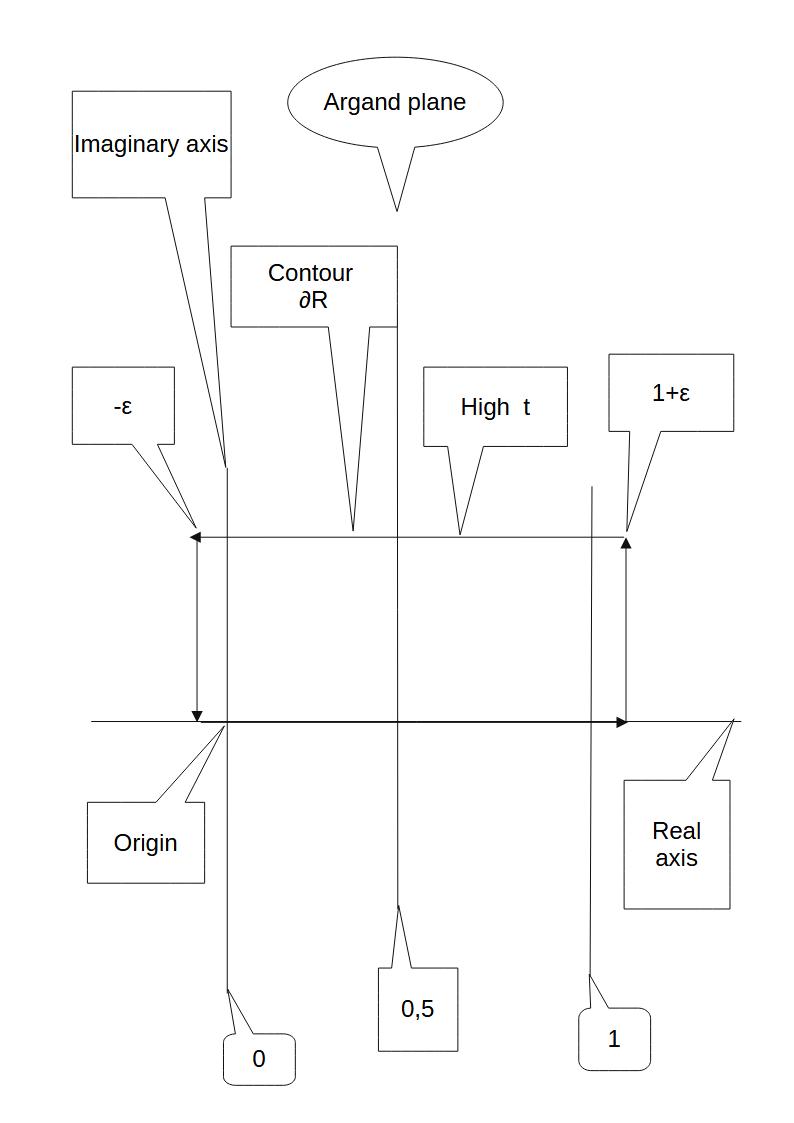}
\caption{Countour $\partial R$}
\label{fig:figura1}
\end{figure}
\\
where it is assumed that t is such that there are no roots $\rho$ on the line $\Im s = t$, and $N(t)$ counts the roots with multiplicities. By symmetry and the fact that $\xi(s)$ is real on the real axis, this can be written as 
\begin{align}
N(t)=\frac{1}{2\pi}\cdot 2 \cdot \Im \left[\int_{C}\frac{\xi'(s)}{\xi(s)}ds \right],\label{3}
\end{align}
where $C$ is the portion of $\partial R$ from $1+\epsilon$ to $\frac{1}{2}+it$.
\begin{figure}[h]
\centering
\includegraphics[width=0.8\textwidth]{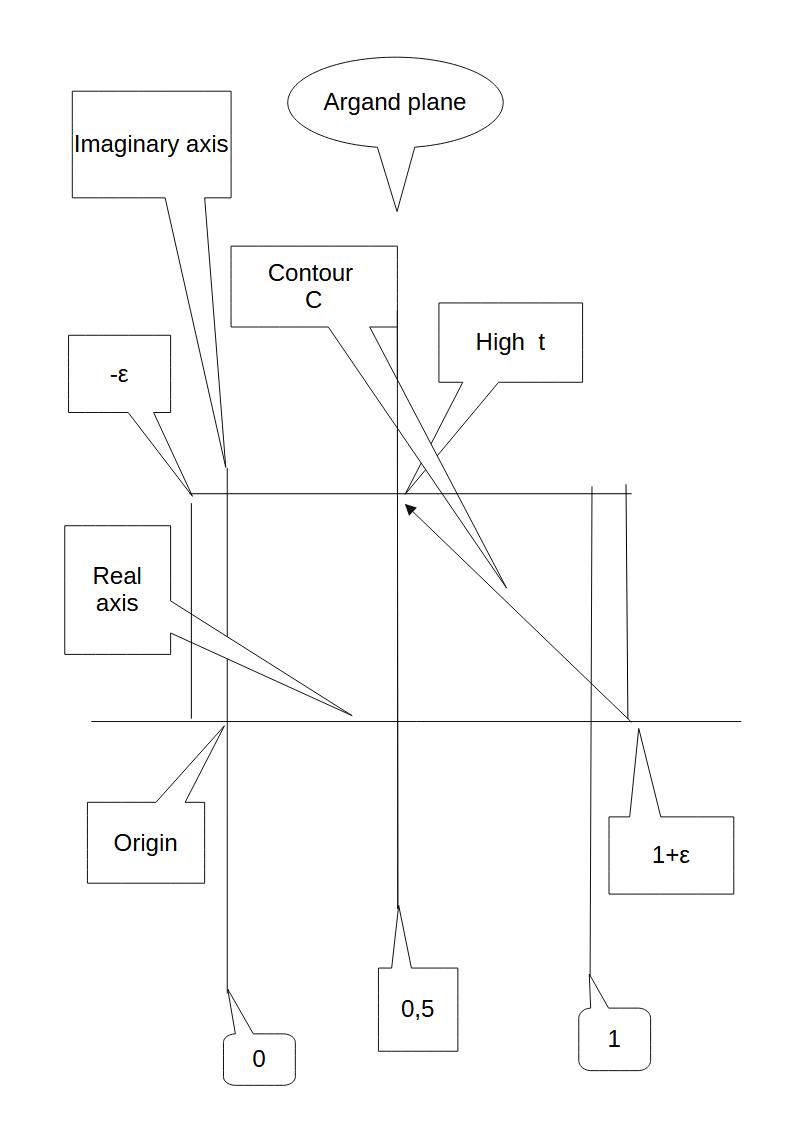}
\caption{Countour $C$}
\label{fig:figura2}
\end{figure}
\\
Using the definition $\xi(s)=\frac{1}{2}s(s-1)\pi^{-s/2}\Gamma(s/2)\zeta(s)$ and the fact that the logarithmic derivative of a product is the sum of the logarithmic derivatives puts this in the form
\begin{align}
&\frac{1}{\pi} \Im \left\{\int_{C} \frac{d}{ds}\ln\left(\pi^{-s/2}\Gamma(s/2)\right)ds \right\}+ \nonumber \\
&\frac{1}{\pi} \Im \left\{\int_{C} \frac{d}{ds}\ln(s(s-1))ds \right\}+ \nonumber \\
&\frac{1}{\pi} \Im \left\{\int_{C}\frac{\zeta'(s)}{\zeta(s)}ds \right\}.\label{4}
\end{align}
The first two terms, being integrals of derivatives, can be evaluated using the fundamental theorem of calculus; the first is  $\pi^{-1}\vartheta(t)$ because it is $\pi^{-1}$ times the imaginary part of $\ln\left( \pi^{-s/2}\Gamma(s/2) \right)$
at $s=\frac{1}{2}+it$ when this $\ln$ is defined to be real on the positive real axis, and the second is 1 because it is $\pi^{-1}$ times the imaginary part of the $\ln$ of $(\frac{1}{2}+it)\cdot (\frac{1}{2}+it-1)=-t^{2}-\frac{1}{4}$ when $\ln(s(s-1))$ is taken to be real for $s>1$. Thus
\begin{align}
N(t)=\frac{1}{\pi}\vartheta(t)+1+\frac{1}{\pi} \Im \int_{C}\frac{\zeta'(s)}{\zeta(s)}ds\label{5}
\end{align}
Backlund observed that if it can be shown that $\Re \zeta(s)$ is never zero on $C$, then this formula suffices to determine N(t) as the integer nearest to $\pi^{-1}\vartheta(t)+1$; this follows simply from noting that if $\Re\zeta(s)$ is never $0$ on $C$, then the curve $\zeta(C)$ never leaves the right halfplane  so that $\ln(\zeta(s))$ is defined al1 along $\zeta(C)$ and gives an antiderivative $\zeta'(s)/\zeta(s)$ on C whose imaginary part lies between $-\pi/2$ and $\pi/2$,which by the fundamental theorem shows that the last term above has absolute value less than $\frac{1}{2}$.

\section{The main result}
We have

\begin{align}
N(t+\delta_{1})-N(t-\delta_{2})=\frac{1}{2 \pi i} \ointctrclockwise_{\partial R_{1}} \frac{\xi'(s)}{\xi(s)}ds
-\frac{1}{2 \pi i} \ointctrclockwise_{\partial R_{2}} \frac{\xi'(s)}{\xi(s)}ds\label{6}
\end{align}
where $\delta_{1}$ and $\delta_{2}$  are infinitesimal quantities that can be chosen, as we will see in the following paragraph, such that no zeros exist on the lines $\Im s = t + \delta_{1}$ and $ \Im s = t - \delta_{2}$. 

Note the graphs below countour $\partial R1.1$ and countour $\partial R1.2$, where circles (the center of the circle is a zero) indicate zeros of the $\zeta(s)$ function, and the distance $d_{min}$ represents the distance from the nearest zero to the line $\Im s = t + \delta_{1}$. Assuming zeros exist on the line $ \Im s = t + \delta_{1}$, with the nearest zero outside the line at distance $d_{min}$, we can choose $\delta_{1}$ to avoid zeros by:

Increasing $\delta_{1}$ (if the nearest zero lies above $ \Im s = t + \delta_{1}$) by a value less than $d_{min}$.

Decreasing $\delta_{1}$ (if the nearest zero lies below $ \Im s = t + \delta_{1}$) by a value less than $d_{min}$.

Similarly, $\delta_{2}$ can be chosen to ensure no zeros occur on the line $\Im s = t - \delta_{2}$.

\begin{figure}[p]
\centering
\includegraphics[width=1\textwidth]{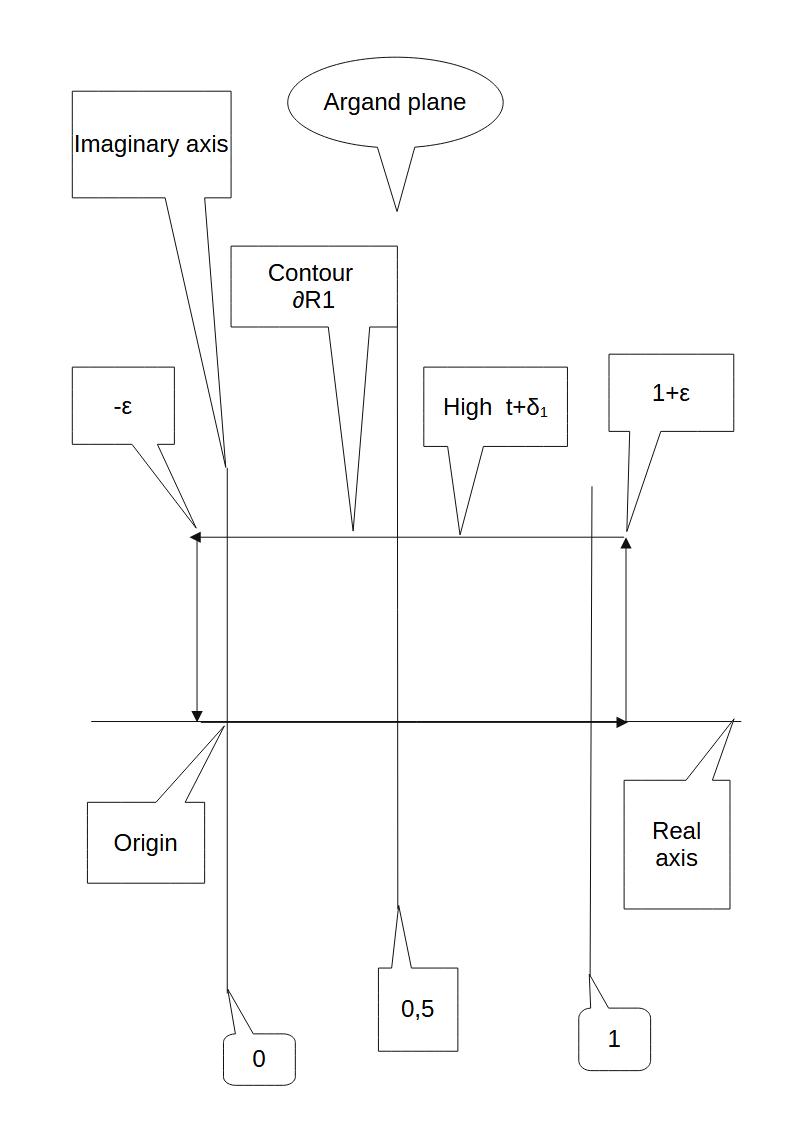}
\caption{Countour $\partial R1$}
\label{fig:figura3}
\end{figure}

\begin{figure}[p]
\centering
\includegraphics[width=1\textwidth]{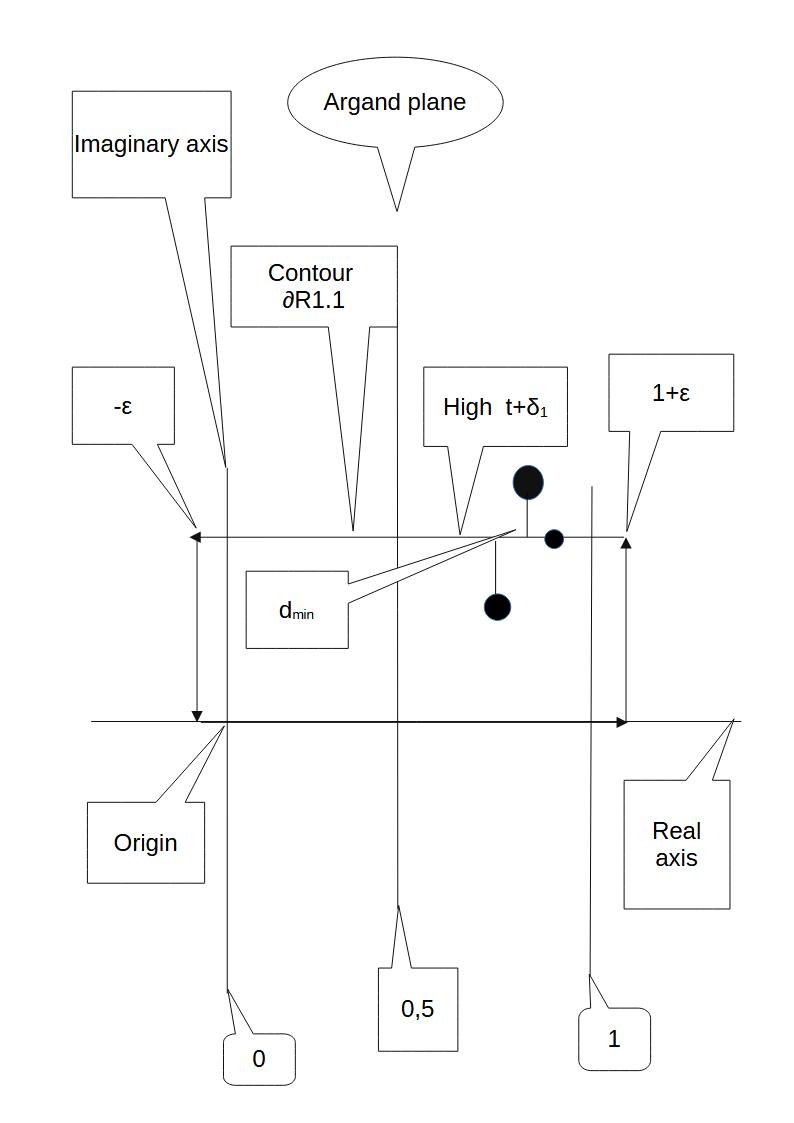}
\caption{Countour $\partial R1.1$}
\label{fig:figura4}
\end{figure}

\begin{figure}[p]
\centering
\includegraphics[width=1\textwidth]{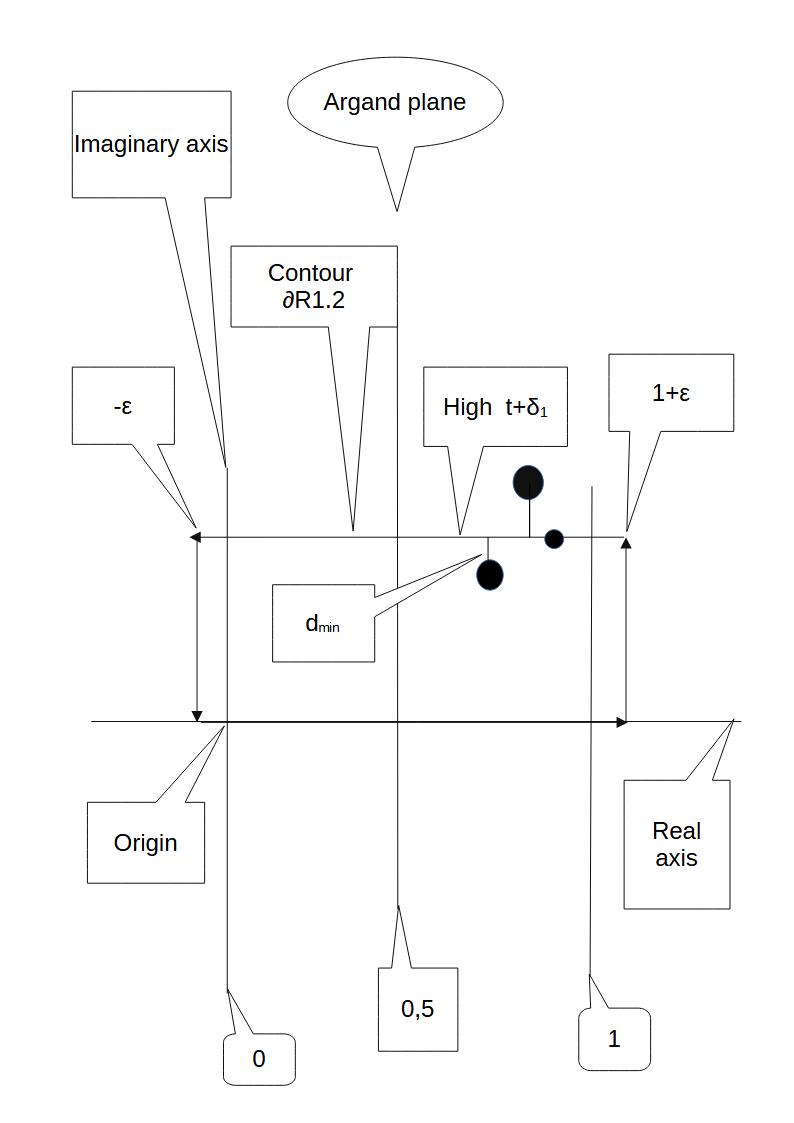}
\caption{Countour $\partial R1.2$}
\label{fig:figura5}
\end{figure}

\begin{figure}[p]
\centering
\includegraphics[width=1\textwidth]{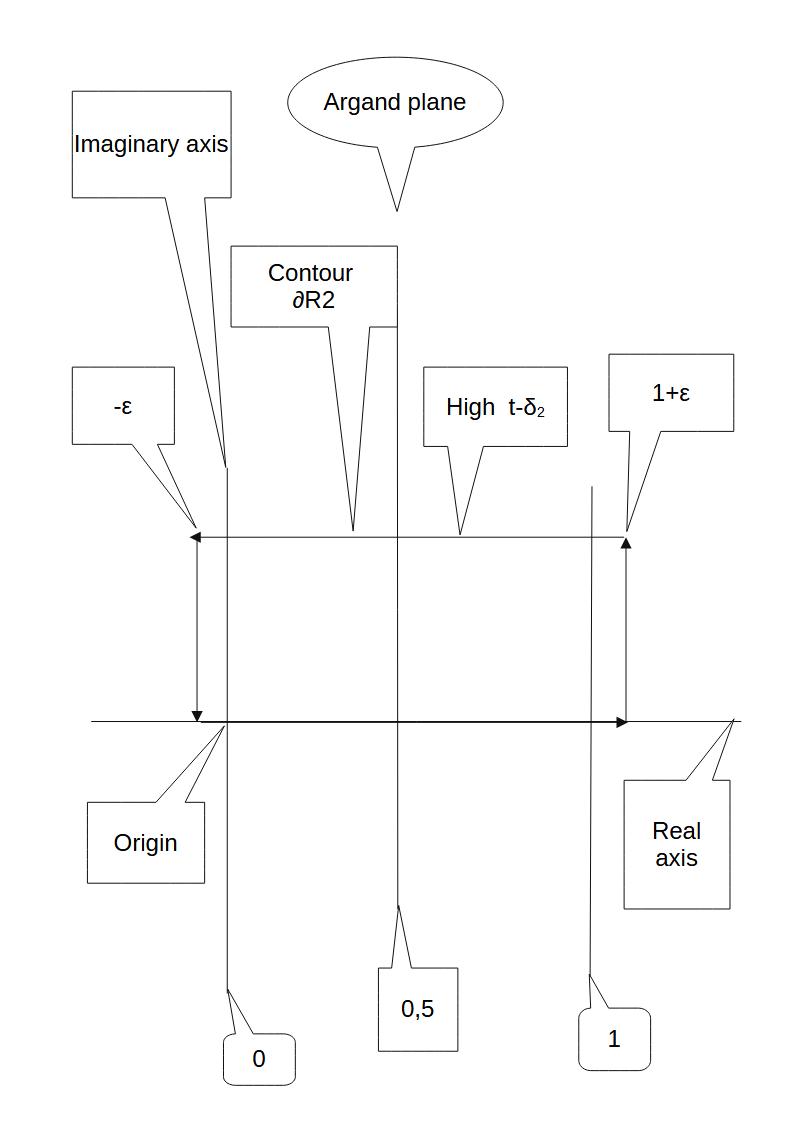}
\caption{Countour $\partial R2$}
\label{fig:figura6}
\end{figure}

\begin{figure}[p]
\centering
\includegraphics[width=0.8\textwidth]{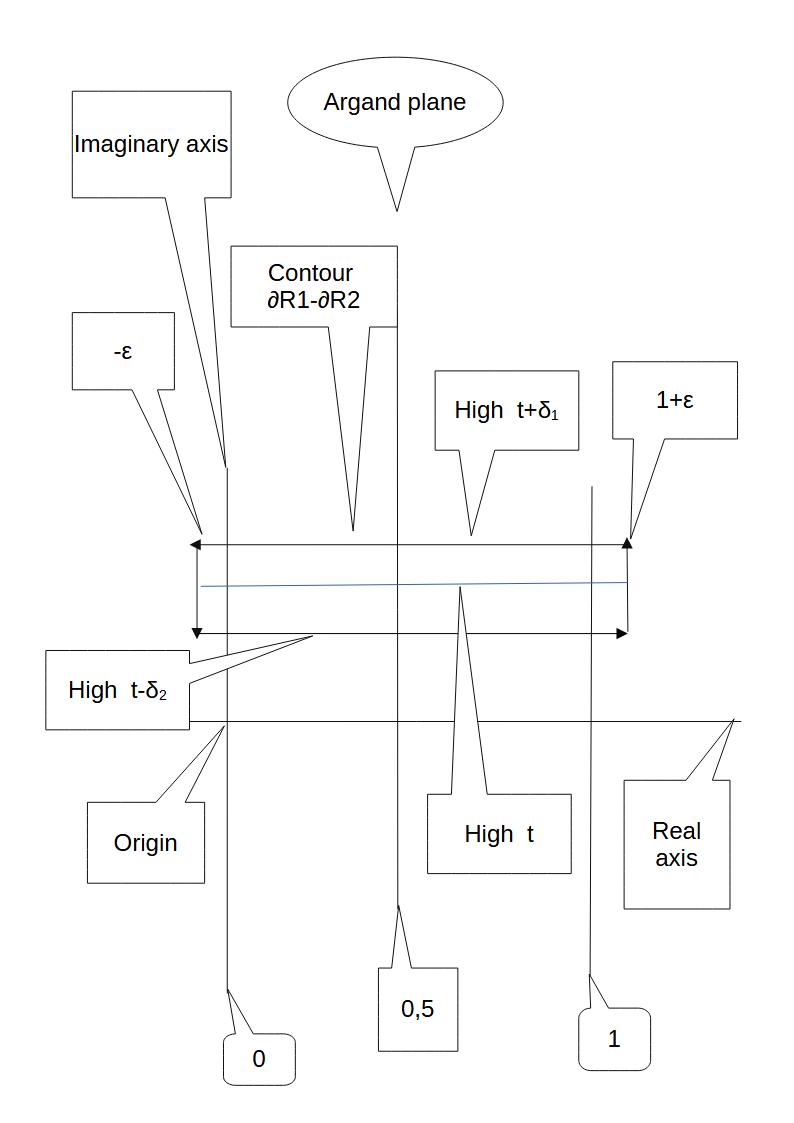}
\caption{Countour $\partial R1-\partial R2$}
\label{fig:figura7}
\end{figure}
\newpage 

\begin{figure}[p]
\centering
\includegraphics[width=1\textwidth]{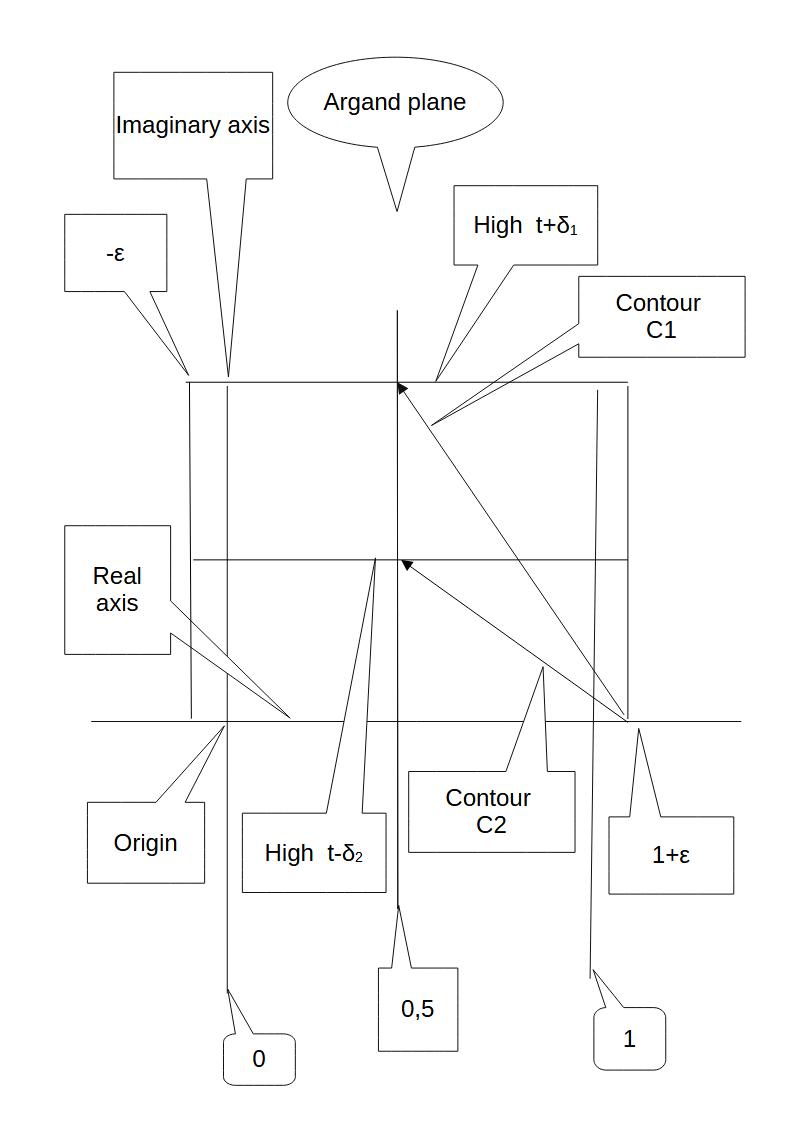}
\caption{Countour $C12$}
\label{fig:figura8}
\end{figure}

\newpage
We have
\begin{align}
N(t+\delta_{1})=\frac{1}{\pi}\vartheta(t+\delta_{1})+1+\frac{1}{\pi} \Im \int_{C_{1}}\frac{\zeta'(s)}{\zeta(s)}ds\label{7}
\end{align}
and
\begin{align}
N(t-\delta_{2})=\frac{1}{\pi}\vartheta(t-\delta_{2})+1+\frac{1}{\pi} \Im \int_{C_{2}}\frac{\zeta'(s)}{\zeta(s)}ds.\label{8}
\end{align}
We can then rewrite equation \ref{6} as
\begin{align}
&N(t+\delta_{1})-N(t-\delta_{2})=\nonumber \\
&\frac{1}{\pi}\vartheta(t+\delta_{1})-\frac{1}{\pi}\vartheta(t-\delta_{2})+ \frac{1}{\pi} \Im \int_{C_{1}}\frac{\zeta'(s)}{\zeta(s)}ds-\frac{1}{\pi} \Im \int_{C_{2}}\frac{\zeta'(s)}{\zeta(s)}ds\label{9}
\end{align}
where $C_{1}$ and $C_{2}$ are the contours shown in the figure 8.

The last two terms on the right-hand side of equation 11 can be estimated. If the real part of the zeta function $\zeta(s)$ does not vanish on lines $C_{1}$ and $C_{2}$, the integrand becomes the logarithmic derivative of zeta, $\frac{\zeta'(s)}{\zeta(s)}$, whose primitive is  $\ln(\zeta(s))$. Thus, the imaginary part is bounded between $-\pi/2$ and $\pi/2$. Hence, we can write

\begin{align}
\frac{1}{\pi} \Im \int_{C_{1}}\frac{\zeta'(s)}{\zeta(s)}ds < \frac{1}{\pi}\cdot \frac{\pi}{2}=\frac{1}{2}\label{10}
\end{align}
\begin{align}
\frac{1}{\pi} \Im \int_{C_{2}}\frac{\zeta'(s)}{\zeta(s)}ds > \frac{1}{\pi}\cdot -\frac{\pi}{2}=-\frac{1}{2}\label{11}
\end{align}
consequently
\begin{align}
\frac{1}{\pi} \Im \int_{C_{1}}\frac{\zeta'(s)}{\zeta(s)}ds-\frac{1}{\pi} \Im \int_{C_{2}}\frac{\zeta'(s)}{\zeta(s)}ds < 1\label{12}
\end{align}
and
\begin{align}
N(t+\delta_{1})-N(t-\delta_{2})< \frac{1}{\pi}\vartheta(t+\delta_{1})-\frac{1}{\pi}\vartheta(t-\delta_{2})+1 < 2.\label{13}
\end{align}
Under the Backlund assumption that $\Re[\zeta(s)] \neq 0$ on contours $C_{1}$ and $C_{2}$, like in equation 7, and choosing $\delta_{1}$ and $\delta_{2}$  infinitesimal therefore $\frac{1}{\pi}\vartheta(t+\delta_{1})-\frac{1}{\pi}\vartheta(t-\delta_{2})<1$, our analysis, if correct, implies the non-existence of two zeros at height t (Figure 7) since possible zeros outside critical line appear in pairs, this count would be at least 2. 

This conclusion extends to all t, indicating zeros lie exclusively on the critical line and have multiplicity one. We will demonstrate in the subsequent section that $\epsilon$ can be selected to preclude zeros on contours $C_{1}$ and $C_{2}$ implying the validity of equation \ref{13}, for all t, i.e. the Riemann Hypothesis.

\section{Selecting epsilon}

By varying epsilon $\epsilon$, we can ensure the real part of zeta  $\Re(\zeta(s))$ doesn't vanish on contours $C_{1}$ and $C_{2}$. 
Contours $C_{1}$ and $C_{2}$ take the form (look figure 8 or 9):

Countour $C_{1}$
\begin{align}
y=-\frac{t+\delta_{1}}{0.5+\epsilon}\cdot (x-0.5)+t+\delta_{1}\label{14}
\end{align}
and  Countour $C_{2}$
\begin{align}
y=-\frac{t-\delta_{2}}{0.5+\epsilon}\cdot (x-0.5)+t-\delta_{2}\label{15}
\end{align}
where the complex argument $s=x+iy$ has real part $x$ and imaginary part $y$.

Suppose that the real part of the zeta function, $\Re \zeta(s)$, vanishes on one or both contours $C_{1}$ and $C_{2}$, at one or more points, and that the real part of the zeta function also vanishes at points near the curves $C_{1}$ and $C_{2}$. Observe the following graphs, where the centers of the black circles represent the zeros of the real part of the zeta function. 
Let $d_{min}$ be defined as the minimum distance between a zero of the real part of the zeta function $\Re \zeta(s)$, located exterior to the contours, and either one or both of the contours $C_{1}$ and $C_{2}$, with the latter case applying when the minimum distance to the contours is equidistant from both.
Two scenarios can occur:
\begin{align}
&\text{1. The nearest zero of the real part is above the curve.} \nonumber\\
&\text{2. The nearest zero of the real part is below the curve.} \nonumber
\end{align}

1 - If the nearest zero is above the line, it is possible to increase $\epsilon$ to shift the line upward, avoiding the external zeros of the real part, since we can select epsilon such that the displacement of all points on the line is less than $d_{min}$, and also those that were initially on the lines (look figure 9). 

2 - If the zero of the real part of zeta is below the line, we can decrease the value of $\epsilon$ to deviate from the external zeros of the lines, since we can again select epsilon such that the displacement of all points on the line is less than $d_{min}$, and also those that were initially on the lines (look figure 10).

\begin{figure}[p]
\centering
\includegraphics[width=1\textwidth]{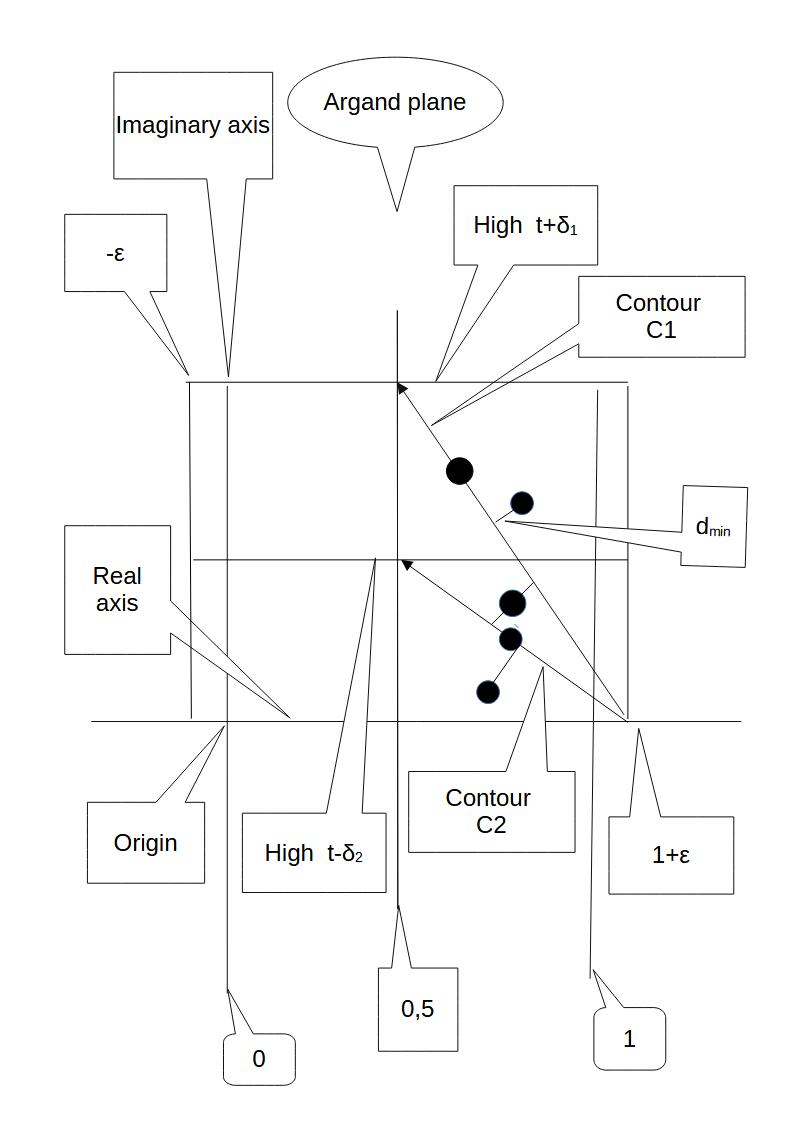}
\caption{Countour $C12.1$}
\label{fig:figura9}
\end{figure}

\begin{figure}[p]
\centering
\includegraphics[width=1\textwidth]{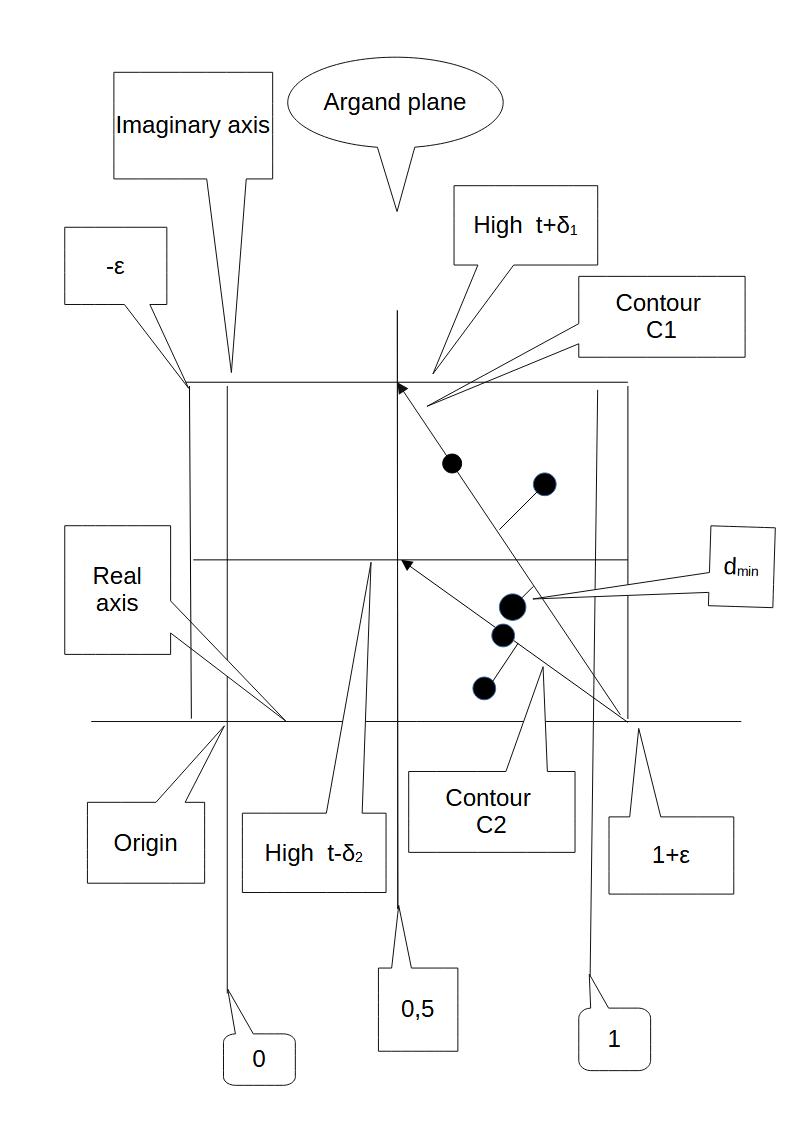}
\caption{Countour $C12.2$}
\label{fig:figura10}
\end{figure}

\begin{figure}[p]
\centering
\includegraphics[width=1\textwidth]{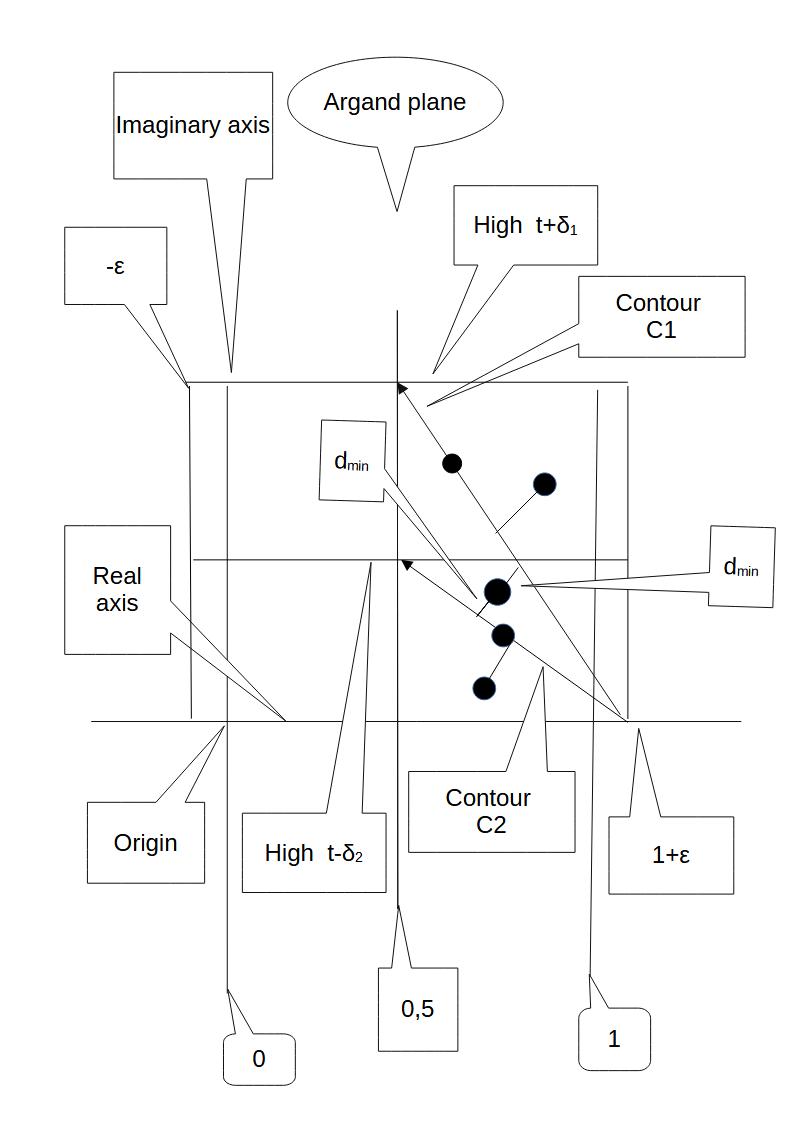}
\caption{Countour $C12.3$}
\label{fig:figura11}
\end{figure}

\begin{figure}[p]
\centering
\includegraphics[width=1\textwidth]{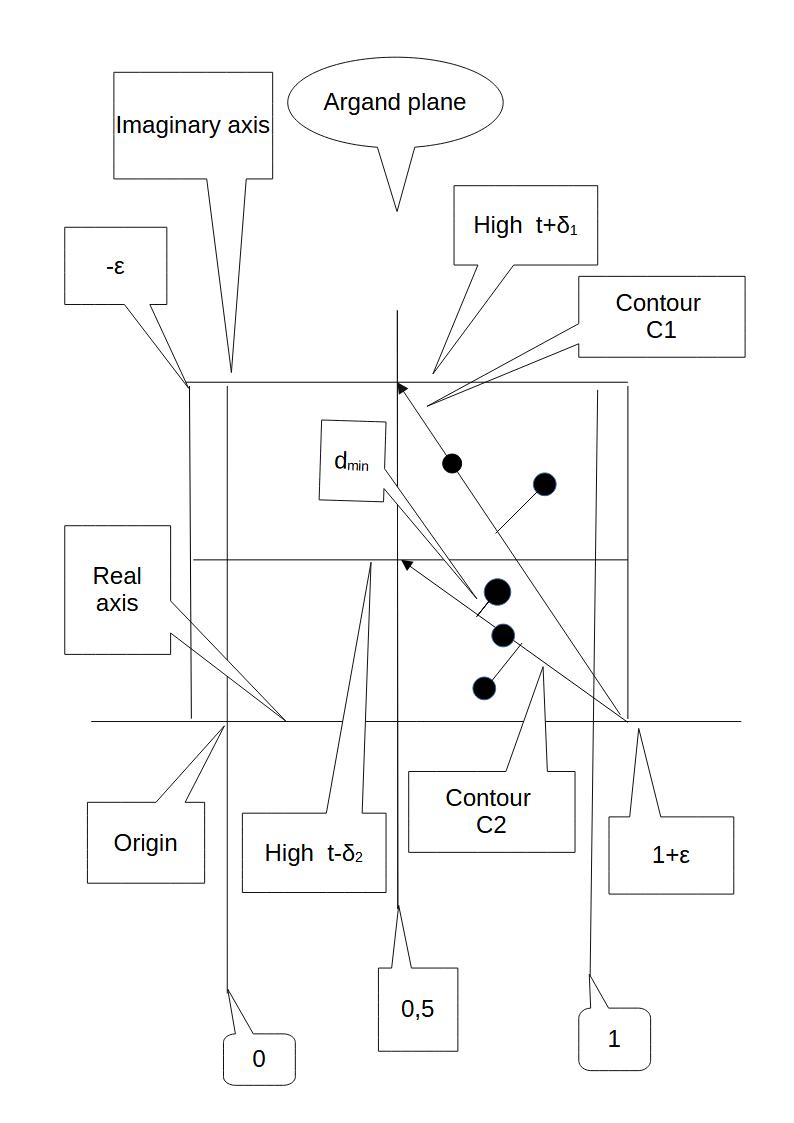}
\caption{Countour $C12.4$}
\label{fig:figura12}
\end{figure}

\begin{figure}[p]
\centering
\includegraphics[width=1\textwidth]{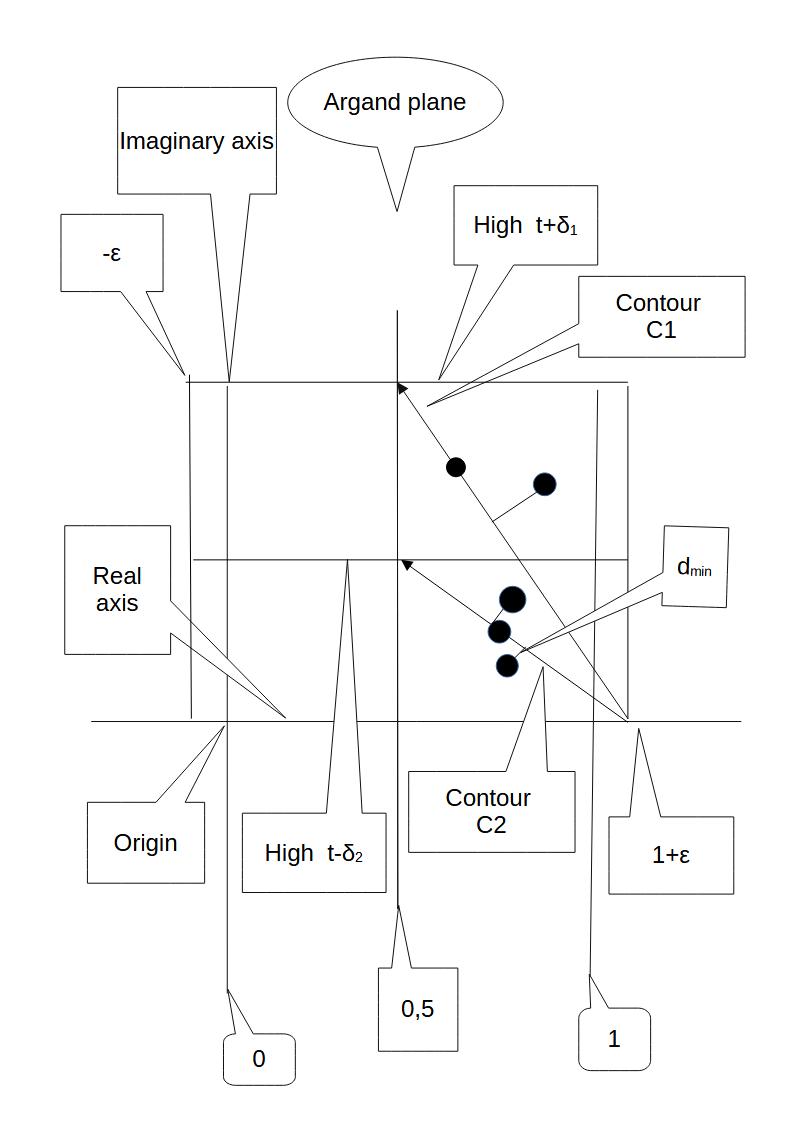}
\caption{Countour $C12.5$}
\label{fig:figura13}
\end{figure}

In the following graphs, we show all possible locations of $d_{min}$:

\begin{align}
&\text{1. $d_{min}$ above $C_{1}$ - Figure \ref{fig:figura9}}\nonumber \\
&\text{2. $d_{min}$ between $C_{1}$ and $C_{2}$, closer to $C_{1}$ - Figure \ref{fig:figura10}} \nonumber \\
&\text{3. $d_{min}$ between $C_{1}$ and $C_{2}$, equidistant from $C_{1}$ and $C_{2}$ - Figure \ref{fig:figura11}}\nonumber \\
&\text{4. $d_{min}$ between $C_{1}$ and $C_{2}$, closer to $C_{2}$ - Figure \ref{fig:figura12}}\nonumber \\
&\text{5. $d_{min}$ below $C_{2}$ - Figure \ref{fig:figura13}} \nonumber
\end{align}

It is demonstrated that with a convenient choice of $\epsilon$, the reasoning of equations (\ref{10}),  (\ref{11}), (\ref{12}) and (\ref{13}) are valid, implying the non-existence of non-trivial zeros outside the critical line, and, if the author has not made a mistake related to the complex analysis, the Riemann Hypothesis is true.

\newpage

\section{Conclusion}
This article presents an attempt to prove the Riemann hypothesis, demonstrating, if the argument is correct, that all non-trivial zeros of the Riemann zeta function satisfy $\sigma=\frac{1}{2}$ . Our approach leverages techniques from complex analysis. The Riemann Hypothesis is equivalent to numerous other problems and has far-reaching implications that are beyond the scope of this article. The present author seeks and desires collaboration and opinions from those interested in solving this important and long-standing problem.

\subsection*{Acknowledgements}
\addcontentsline{toc}{subsection}{Dedication}
\textit{I dedicate this work to my beloved parents, Rodolpho Antonio de Rezende and Edna Vieira Rocha de Rezende (Posthumously), and to my equally dear siblings, Gustavo Rocha de Rezende and Gisella Rocha de Rezende.}

\end{document}